\documentclass[12pt,
amscd]{amsart}
\usepackage{latexsym,amscd,amssymb}  
\theoremstyle{plain}
  \newtheorem{thm}{Theorem}[section]
  \newtheorem{prop}[thm]{Proposition}
  \newtheorem{lem}[thm]{Lemma}
  \newtheorem{cor}[thm]{Corollary}
\theoremstyle{definition}
  \newtheorem{dfn}[thm]{Definition}
  \newtheorem{exmp}[thm]{Example}
\theoremstyle{remark}
  \newtheorem{rem}[thm]{Remark}
\def\ba{{\bf a}}
\def\bb{{\bf b}}

\def\C{C^\bullet}
\def\D{\omega^\bullet}
\def\E{{}^*E}
\def\bT{{\mathbf T}}
\def\I{I^\bullet}
\def\J{J^\bullet}
\def\M{M^\bullet}
\def\N{N^\bullet}
\def\P{P^\bullet}
\def\Pro{\bf P}
\def\Q{Q^\bullet}
\def\bx{x}
\def\SS{\operatorname{SS}}
\def\lk{\operatorname{lk}}
\def\Syz{\operatorname{Syz}}

\def\Rep{\operatorname{\bf Rep}}
\def\rH{\tilde{H}}
\def\op{{\rm op}}
\def\gr{\operatorname{gr}_\m}
\def\grF{\operatorname{gr}}
\def\m{{\mathfrak m}}
\def\rad{{\mathfrak r}}
\def\NN{{\mathbb N}}
\def\Sq{{\bf Sq}}
\def\SqE{{{\bf Sq}_E}}
\def\SqS{{\bf Sq}_S}
\def\Str{{\bf Str}}
\def\ZZ{{\mathbb Z}}
\def\RR{{\mathbb R}}

\def\cD{{\bf D}}
\def\cA{{\bf A}}

\def\cS{{\mathcal S}}
\def\cE{{\mathcal E}}
\def\cL{{\mathcal L}}
\def\cR{{\mathcal R}}
\def\bP{{\mathcal F}}

\def\Kom{\operatorname{Com}}
\def\Nc{{\mathcal N}}
\def\Zc{{\mathcal Z}}
\def\b1{{\mathbf 1}}
\def\Hom{\operatorname{Hom}}
\def\can{\omega_S}
\def\Ext{\operatorname{Ext}}

\def\chara{\operatorname{char}}
\def\cmpl{{\sf c}}

\def\supp{\operatorname{supp}_+}

\def\<{{\langle}}
\def\>{{\rangle}}
\def\md{{\bf mod}_\Lambda}
\def\mdZ{{\bf gr.mod}_\Lambda}
\def\Amd{{\bf mod}_A}
\def\mdR{{\bf gr.mod}_R}
\def\mmR{{\bf gr.mod}_{R^!}}
\def\mmS{{\bf mod}_S}
\def\mmE{{\bf mod}_E}
\def\MMZn{{{}^*{\bf Mod}}}

\def\mmEn{{}^*{\bf mod}_E}
\def\mmSn{{}^*{\bf mod}_S}
\def\MM{{\bf Mod}}

\def\too{\longrightarrow}
\def\Mus{Musta\c{t}\v{a}}

\def\rat{{\rm rat}}
\def\Id{\operatorname{Id}}
\def\injdim{\operatorname{inj.dim}}
\def\pd{\operatorname{proj.dim}}
\begin{document}
\title[Derived Category of Squarefree Modules]
{Derived Category of Squarefree Modules and Local Cohomology with Monomial 
Ideal Support}
\author{Kohji Yanagawa}
\address{Department of Mathematics, 
Graduate School of Science, Osaka University, Toyonaka, Osaka 
560-0043, Japan}
\email{yanagawa@math.sci.osaka-u.ac.jp}
\keywords{Stanley-Reisner ring, local duality, Alexander duality,  
local cohomology, Bernstein-Gel'fand-Gel'fand correspondence, Koszul duality}  
\subjclass{Primary 13D25; Secondary 13D02, 13D45, 13F55,  
13N10, 18E30}
\maketitle
\begin{abstract}
A {\it squarefree module} over a polynomial ring $S = k[x_1, \ldots, x_n]$ is 
a generalization of a Stanley-Reisner ring, and allows us to apply homological 
methods to the study of monomial ideals more systematically. 

The category $\Sq$ of squarefree modules 
is equivalent to the category of finitely generated 
left $\Lambda$-modules, where $\Lambda$ is the incidence algebra 
of the Boolean lattice $2^{\{1, \ldots, n\}}$. 
The derived category $D^b(\Sq)$ has two duality functors 
$\cD$ and $\cA$. The functor $\cD$ is a common one 
with $H^i(\cD(\M)) = \Ext_S^{n+i}(\M,\can)$, 
while the {\it Alexander duality functor} 
$\cA$ is rather combinatorial. We have a strange relation 
$\cD \circ \cA \circ \cD \circ \cA \circ \cD \circ \cA 
\cong \bT^{2n}$, where $\bT$ is the translation functor. 
The functors $\cA \circ \cD$ and $\cD \circ \cA$ 
give a non-trivial autoequivalence of $D^b(\Sq)$. 
This equivalence corresponds to the Koszul duality for $\Lambda$,  
which is a Koszul algebra with $\Lambda^! \cong \Lambda$.  
Our $\cD$ and $\cA$ are also related to 
the Bernstein-Gel'fand-Gel'fand correspondence. 

The local cohomology  $H_{I_\Delta}^i(S)$ at a Stanley-Reisner 
ideal $I_\Delta$ 
can be constructed from the squarefree module $\Ext_S^i(S/I_\Delta, \can)$. 
We see that Hochster's formula on the $\ZZ^n$-graded Hilbert function of 
$H_\m^i(S/I_\Delta)$ is also a formula on the characteristic cycle of 
$H_{I_\Delta}^{n-i}(S)$ as a module over the Weyl algebra
$A = k\<x_1, \ldots, x_n,$ $\partial_1, \ldots, \partial_n \>$ 
(if $\chara(k)=0$).
\end{abstract}

\section{Introduction}
The Stanley-Reisner ring of an abstract simplicial complex $\Delta \subset 
2^{\{1, \ldots, n\}}$, which is the quotient of a polynomial ring 
$S = k[x_1, \ldots, x_n]$ by the squarefree monomial ideal $I_\Delta$, 
is a central concept of combinatorial commutative algebra (see
\cite{BH,St}). In \cite{Y}, the author defined a {\it squarefree} 
$\NN^n$-graded $S$-module.  A Stanley-Reisner ring $S/I_{\Delta}$, 
its syzygy module $\Syz_i(S/I_\Delta)$, the canonical module $\can$, and
$\Ext^i_S(S/I_\Delta, \omega_S)$ are always squarefree. Using this notion,  
we can treat Stanley-Reisner rings and related objects in a categorical way, 
see \cite{EPY,Mil,RWY,R0,R,Y1,Y2,Y3}.  In the present paper, 
we will study the derived category of squarefree modules. 

Let $\SqS$ (or simply $\Sq$)
be the category of squarefree $S$-modules and 
their degree preserving maps. Then $\Sq$ is equivalent to the category of 
finitely generated left $\Lambda$-modules, where $\Lambda$ is the incidence 
algebra of the Boolean lattice $2^{\{1, \ldots, n\}}$. 

Let $D^b(\Sq)$ be the derived category of bounded complexes in $\Sq$. 
We have contravariant functors $\cD$ and $\cA$ from $D^b(\Sq)$ to 
itself satisfying $\cD^2 \cong  \cA^2 \cong \Id_{D^b(\Sq)}$. 
If $M$ is a squarefree module, so is  $\Ext_S^i(M, \can)$. 
Moreover, for $\M \in D^b(\Sq)$, we can define $\cD(\M) \in D^b(\Sq)$ 
with $H^i(\cD(\M)) = \Ext_S^{n+i}(\M,\can)$ 
in a natural way, see Proposition~\ref{Ext}. On the other hand, 
extending an idea of Eagon-Reiner~\cite{ER},  Miller~\cite{Mil} and 
R\"omer~\cite{R0} constructed the {\it Alexander duality functor}  $\cA$ on 
$\Sq$. Since $\cA$ is exact, we can regard it as a duality 
functor on $D^b(\Sq)$. 

Using $D^b(\Sq)$, we can get simple and systematic proofs of many results 
in \cite{Mil,R0,R,Y,Y1}.  
Moreover, we  prove a strange natural equivalence  
$$\cD \circ \cA \circ \cD \circ \cA \circ \cD \circ \cA \cong \bT^{2n},$$ 
where $\bT$ is the translation functor on $D^b(\Sq)$. 

Let $E =  \bigwedge S_1^*$ be the exterior algebra. 
A squarefree module over $E$, which was defined by R\"omer \cite{R0}, 
is also a natural concept. The category $\SqE$ of squarefree $E$-modules 
is equivalent to $\SqS$ in a natural way. 
A famous theorem of Bernstein-Gel'fand-Gel'fand~\cite{BGG} states that 
the bounded derived category of finitely generated $\ZZ$-graded $S$-modules 
is equivalent to the bounded derived category of finitely generated 
$\ZZ$-graded left $E$-modules. 
The functors defining this equivalence preserve the squarefreeness, 
and coincide with $\cA \circ \cD$ and $\cD \circ \cA$ in the squarefree case 
under the equivalence $\SqS \cong \SqE$. 
We have another relation to Koszul duality. 
The incidence algebra $\Lambda$ of $2^{\{1, \ldots, n\}}$ 
is a Koszul algebra whose quadratic dual $\Lambda^!$ is isomorphic to 
$\Lambda$ itself. The functors $\cA \circ \cD$ and $\cD \circ \cA$ 
give a non-trivial autoequivalence of $D^b(\Sq)$. 
This equivalence corresponds to the Koszul duality 
$D^b(\md) \cong D^b({\bf mod}_{\Lambda^!})$. 

In the last section, under the assumption that $\chara(k) = 0$,  
we study modules over the Weyl algebra 
$k\<x_1, \ldots, x_n, \partial_1, \ldots, \partial_n \>$ 
associated to squarefree modules (e.g., the local cohomology module 
$H_{I_\Delta}^i(S)$). Especially, 
we give the formula for their characteristic cycles. 

\medskip

After I received the referee's report for the first version, 
I widely revised the paper. Among other things, 
Proposition~\ref{self Koszul dual} is a new result of the second version 
which was submitted in September 2001. 
The present version is the fourth one, in which some proofs and expositions 
are revised.

\section{Preliminaries}
Let $S = k[x_1, \ldots, x_n]$ be a polynomial ring over a field $k$. 
Consider an $\NN^n$-grading $S = \bigoplus_{\ba \in \NN^n} S_\ba = 
\bigoplus_{\ba \in \NN^n} k \, \bx^\ba$, where $\bx^\ba = \prod_{i=1}^n 
x_i^{a_i}$ is the monomial with the exponent $\ba = (a_1, \ldots, a_n)$. 
We denote the graded maximal ideal $(x_1, \ldots, x_n)$ by $\m$. 
For a $\ZZ^n$-graded module $M$ and $\ba \in \ZZ^n$, $M_\ba$ 
means the degree $\ba$ component of $M$, and $M(\ba)$ denotes the 
shifted module with $M(\ba)_\bb = M_{\ba+\bb}$. 
We denote the category of $S$-modules by $\MM$, 
and the category of $\ZZ^n$-graded $S$-modules by $\MMZn$.  
Here a morphism $f$ in $\MMZn$ is an $S$-homomorphism 
$f: M \to N$ with $f(M_\ba) \subset N_\ba$ for all $\ba \in \ZZ^n$.  
See \cite{GW} for information on $\MMZn$. 

For $M,N \in \MMZn$ and $\ba \in \ZZ^n$, set  
${}^*\Hom_S(M, N)_\ba := \Hom_\MMZn(M,N(\ba))$.  
Then $${}^*\Hom_S(M,N) := \bigoplus_{\ba \in \ZZ^n}{}^*\Hom_S(M, N)_\ba$$ 
has a natural $\ZZ^n$-graded $S$-module structure. 
If $M$ is finitely generated, then  ${}^*\Hom_S(M,N)$ is isomorphic to 
the usual $\Hom_S(M,N)$ as the underlying $S$-module. 
Thus, we simply denote ${}^*\Hom_S(M,N)$ by $\Hom_S(M,N)$ in this case. 
In the same situation, $\Ext_S^i(M,N)$ also has a $\ZZ^n$-grading with 
$\Ext_S^i(M,N)_\ba = \Ext^i_\MMZn(M,N(\ba))$. 

For $\ba \in \ZZ^n$, set $\supp (\ba) := \{i  \mid a_i > 0\} \subset [n] 
:= \{1, \ldots, n \}$. 
We say  $\ba \in \ZZ^n$ is {\it squarefree} if $a_i= 0,1$ for all $i \in [n]$. 
When $\ba \in \ZZ^n$ is squarefree, we sometimes identify $\ba$ with  
$\supp (\ba)$. Let $\Delta \subset 2^{[n]}$ be a simplicial complex   
(i.e., $\Delta \ne \emptyset$, and 
$F \in \Delta$ and $G \subset F$ imply $G \in \Delta$). 
The {\it Stanley-Reisner ideal} of $\Delta$ is the squarefree monomial ideal 
$I_{\Delta} := (\bx^F \, | \, F \not \in \Delta)$ of $S$. Any squarefree 
monomial ideal is the Stanley-Reisner ideal $I_\Delta$ for some $\Delta$. 
We say $S/I_\Delta$ is the {\it Stanley-Reisner ring} of $\Delta$. 

\begin{dfn}[\cite{Y}]
We say a $\ZZ^n$-graded $S$-module $M$ is {\it squarefree}, if 
the following conditions are satisfied. 
\begin{itemize}
\item[(a)] $M$ is $\NN^n$-graded (i.e., $M_\ba =0$ if $\ba \not \in 
\NN^n$), and $\dim_k M_\ba  < \infty$ for all $\ba \in \NN^n$. 
\item[(b)] The multiplication map 
$M_{\ba} \ni y \mapsto \bx^\bb y \in M_{\ba +\bb}$ is 
bijective for all $\ba, \bb \in \NN^n$ with $\supp (\ba+\bb) = \supp (\ba)$. 
\end{itemize}
\end{dfn}

A squarefree module $M$ is generated by its squarefree part 
$\bigcup_{F \subset [n]} M_F$. Thus it is finitely generated.  
For a simplicial complex $\Delta \subset 2^{[n]}$, 
$I_\Delta$ and  $S/I_{\Delta}$ are squarefree modules. 
A free module $S(-F)$, $F \subset [n]$, is also squarefree.  In particular, 
the $\ZZ^n$-graded canonical module $\can = S(-\b1)$ of $S$ is squarefree, 
where $\b1 = (1, \ldots, 1)$. 

Denote by $\Sq_S$ (or simply $\Sq$)
the full subcategory of $\MMZn$ consisting 
of all the squarefree modules.  In $\MMZn$, $\Sq$ is closed under 
kernels, cokernels and  extensions (\cite[Lemma~2.3]{Y}). 
For the study of $\Sq$, the incidence algebra 
of a finite partially ordered set ({\it poset}, for short) is very useful,  
as shown in \cite{RWY,Y3}.  In Section 4 of the present paper, we will use 
further properties of the incidence algebra (of a Boolean lattice). 
So we now recall basic properties of an incidence algebra 
for the reader's convenience. 
See \cite[\S III.1]{ARS} for undefined terminology.

Let $P$ be a finite poset. The incidence algebra $\Lambda = I(P,k)$ 
of $P$ over $k$ is the $k$-vector space with a basis 
$\{e_{x, \, y} \mid \text{$x, y \in P$ with $x \geq y$} \}$. 
The $k$-bilinear multiplication defined  by $e_{x, \, y} \, e_{z, \, w} = 
\delta_{y, \, z} \, e_{x, \, w}$ makes $\Lambda$ 
a finite dimensional associative $k$-algebra. (The usual definition is 
the opposite ring of our $\Lambda$. But we use the above definition for the 
convenience in a later section.) Set $e_x := e_{x, \, x}$. 
Then $1 = \sum_{x \in P} e_x$ and $e_x \, e_y = \delta_{x,y} \, e_x$. We have 
$\Lambda \cong \bigoplus_{x \in P} \Lambda e_x$ as a left $\Lambda$-module, 
and each $\Lambda e_x$ is indecomposable. 

An incidence algebra $\Lambda$ is the algebra associated 
with a {\it quiver with relations}. For 
a poset $P$, we consider the quiver $\Gamma = \{ \Gamma_0, \Gamma_1\}$ 
with $\Gamma_0 = P$ and 
$$\Gamma_1 = \{ \, x \, \cdot  \longleftarrow  \cdot \, y \mid 
\text{$x , y \in P$, $x > y$, but there is {\it no} 
$z \in P$ with $x > z> y$} \, \}.$$
So $\Gamma$ is (essentially) the Hasse diagram of $P$. 
Set 
$$\rho := \{ \, p_1 - p_2 \mid \text{$p_1$ and $p_2$ are paths of $\Gamma$ 
with $s(p_1) = s(p_2)$ and $e(p_1) = e(p_2)$} \, \},$$
where $s(p_i)$ and $e(p_i)$ represent the initial vertex and the final vertex 
of $p_i$ respectively. Let $k(\Gamma, \rho)$ be the algebra associated 
with $(\Gamma, \rho)$. Then we have an isomorphism $\psi : k(\Gamma, \rho) 
\stackrel{\cong}{\to} \Lambda$. Here, if $[p]$ is the residue class  
containing a path $p$ of $\Gamma$, we have $\psi([p]) = e_{x,\, y}$, 
where $x = e(p)$ and $y = s(p)$. 

Denote the category of finitely generated left $\Lambda$-modules by $\md$. 
If $N \in \md$, we have $N = \bigoplus_{x \in P} N_x$ as a 
$k$-vector space, where $N_x := e_x N$.  Note that $e_{x, \, y} \, 
N_y \subset N_x$ and $e_{x,\, y} \, N_z  = 0$ for $y \ne z$. 
If $f: N \to N'$ is a morphism in $\md$, then $f(N_x) \subset N'_x$.  
Under the isomorphism $\Lambda \cong k(\Gamma, \rho)$, 
$\md$ is equivalent to  the category $\Rep(\Gamma, \rho)$ of representations 
of $(\Gamma, \rho)$ by \cite[III, Proposition~1.7]{ARS}. 
If $(V, f) \in \Rep(\Gamma, \rho)$ corresponds to 
$N \in \md$, then $N_x = V(x)$ for $x \in P$. 
We have explicit descriptions of simple objects, indecomposable projectives, 
and indecomposable injectives in 
$\md \cong \Rep(\Gamma, \rho)$, see \cite[\S III. 1]{ARS}. 

Let $2^{[n]}$ be the Boolean lattice (i.e., we regard the power set 
$2^{[n]}$ of $[n]$ as a poset by inclusions), and $\Lambda= I(2^{[n]}, k)$ 
its incidence algebra. 
For $M \in \Sq$, set $\Phi(M) := N = \bigoplus_{F \subset [n]} N_F$ to be a 
$k$-vector space with $M_F \cong N_F$. 
Then $N$ has a left $\Lambda$-module structure such that the 
multiplication map $N_F \ni y \mapsto e_{G, \, F} \, y \in N_{G}$ 
for $G \supset F$ is induced by $M_F \ni y \mapsto \bx^{(G \setminus F)} 
y  \in M_G$. It is easy to see that $\Phi$ gives a covariant functor 
$\Sq \to \md$. Recall that $\Lambda \cong k(\Gamma, \rho)$, 
where $\Gamma$ is a quiver whose set of vertices is 
$2^{[n]}$, and $\md \cong \Rep(\Gamma, \rho)$. If $M$ is a squarefree module, 
$\Phi(M)$ corresponds to the representation $(V, f) \in \Rep(\Gamma, \rho)$ 
with $V(F) = M_F$ and $f_{F \cup \{ i \},\, F}: 
V(F)=M_F \ni y \mapsto x_i y \in M_{F \cup \{ i \}} = V(F \cup \{i \})$ 
for $F \subset [n]$ and $i \in [n] \setminus F$. 
In \cite{Y2}, the author used {\it sheaves on a poset} 
to understand squarefree modules. But this notion is equivalent to 
that of representations of $(\Gamma, \rho)$ in our context.

\begin{prop}[{\cite{Y2,Y3}}]\label{sq=sh}
Let $\Lambda= I(2^{[n]}, k)$ be the incidence algebra. 
The functor $\Phi$ constructed above gives 
an equivalence $\Sq \cong \md$. 
\end{prop}

For a subset $F \subset [n]$, $P_F$ denotes the monomial prime ideal 
$(x_i \mid i \not \in F)$ of $S$. The next result follows from 
Proposition~\ref{sq=sh} and \cite[\S III. 1]{ARS}. 

\begin{cor}[\cite{Y1}]\label{proj inj}
$\Sq$ is an abelian category, and has enough projectives and 
injectives. An indecomposable projective (resp. injective) object in $\Sq$ 
is isomorphic to $S(-F)$ (resp. $S/P_F$) for some $F \subset [n]$. 
For any squarefree module $M$, both $\pd_\Sq M$ and $\injdim_\Sq M$ are 
at most $n$. 
\end{cor}

Many invariants of squarefree modules are naturally 
described in terms of $\Lambda$. For example, 
if $M$ is a squarefree module with $N := \Phi(M)$, 
$\dim_S M = \max \{ \,  |F| \mid N_F \ne 0 \, \} 
= n - \min \{ \, i \mid \Ext_\Lambda^i(N, \Lambda) \ne 0 \, \}$ and 
$\pd_S M = \pd_\Lambda N = \max 
\{ \, i \mid \Ext_\Lambda^i(N, \Lambda) \ne 0 \, \}$. 
See Remark~\ref{Ext_A} below for information on 
$\Ext_\Lambda^i(N, \Lambda)$.  

We also remark that $\Sq$ admits the Jordan-H\"older 
theorem and the Krull-Schmidt theorem and a simple object in $\Sq$ 
(i.e., a non-zero squarefree module without non-trivial squarefree submodule)
is isomorphic to $(S/P_F)(-F)$ for some $F$. 

\begin{dfn}[\cite{Y1}]\label{straight}
A $\ZZ^n$-graded $S$-module $M = \bigoplus_{\ba \in \ZZ^n} M_\ba$ is called 
{\it straight}, if the following two conditions are satisfied. 
\begin{itemize}
\item[(a)] $\dim_k M_{\ba} < \infty$ for all $\ba \in \ZZ^n$. 
\item[(b)] The multiplication map 
$M_{\ba} \ni y \mapsto \bx^\bb y \in M_{\ba + \bb}$ is bijective 
for all $\ba \in \ZZ^n$ and $\bb \in \NN^n$ with $\supp 
(\ba + \bb) = \supp (\ba)$. 
\end{itemize}
\end{dfn}

For a $\ZZ^n$-graded $S$-module $M = \bigoplus_{\ba \in \ZZ^n} 
M_\ba$,  we call the submodule $\bigoplus_{\ba \in \NN^n} 
M_\ba$ the {\it $\NN^n$-graded part} of $M$, and denote it by $\Nc(M)$. 
If $M$ is straight then $\Nc(M)$ is squarefree. 
Conversely, for any squarefree module $N$, there is a unique 
(up to isomorphism) straight module $\Zc(N)$ whose $\NN^n$-graded 
part is isomorphic to $N$. For example, $\Zc(S/P_F) \cong \E(S/P_F)$, 
where $\E(S/P_F)$ is the injective envelope of $S/P_F$ in $\MMZn$. 
Denote by $\Str_S$ (or simply $\Str$) the full subcategory of $\MMZn$ 
consisting of all the straight $S$-modules.  

\begin{prop}[{\cite[Proposition~2.7]{Y1}}] 
The functors $\Nc: \Str \to \Sq$ and $\Zc: \Sq \to \Str$ 
give an equivalence $\Sq \cong \Str$.  
\end{prop}

\medskip

Let $\Kom^b(\Sq)$ be the category of bounded cochain complexes of 
squarefree modules, and $D^b(\Sq)$ the bounded derived category of $\Sq$. 
A squarefree module $M$ can be regarded as a complex 
$\cdots \to 0 \to M \to 0 \to \cdots $ 
with $M$ at the $0^{\rm th}$ place. 
For a complex $\M$ and an integer $p$, let $\M[p]$ be the $p^{\rm th}$ 
{\it translation} of $\M$. That is, $\M[p]$ is a complex with 
$M^i[p] = M^{i+p}$ and $d_{M[p]} = (-1)^pd_M$. 

A complex $\M \in \Kom^b(\Sq)$ has 
a projective resolution $\P \in \Kom^b(\Sq)$. That is,  
there is a quasi-isomorphism $\P \to \M$ and each $P^i$ is  
projective in $\Sq$. We say $\P$ is {\it minimal} if 
$d_P (P^{i-1}) \subset \m P^i$ for all $i$. 
A minimal projective resolution of $\M \in \Kom^b(\Sq)$ in $\Sq$ 
is a $\ZZ^n$-graded minimal $S$-free resolution of $\M$. 
Under the same notation as Proposition~\ref{sq=sh},  
a projective resolution $\P \in \Kom^b(\Sq)$ is minimal 
if and only if so is $\Q :=\Phi(\P) \in \Kom^b(\md)$, that is, 
$d_Q (Q^{i-1}) \subset \rad Q^i$ for all $i$. 
Here $\rad = \<e_{F, \, G} \mid  F \supsetneq G \>$ 
is the Jacobson radical of $\Lambda$. 
Hence every $\M \in \Kom^b(\Sq)$ has a unique minimal 
projective resolution, and any projective resolution is a direct 
sum of a minimal one and an exact complex.   
Let $\P$ be a minimal projective resolution of $\M \in \Kom^b(\Sq)$. 
We define $\beta_i(F, \M) \in \NN$ so that 
$$P^{-i} \cong \bigoplus_{F \subset [n]}S(-F)^{\beta_i(F, \M)}.$$

Similarly, every $\M \in \Kom^b(\Sq)$ has an injective resolution 
$\I \in \Kom^b(\Sq)$. That is,  there is a quasi-isomorphism $\M \to \I$ 
and each $I^i$ is injective in $\Sq$. 
We say $\I$ is {\it minimal} if $I^i$ is a *essential extension 
of $\ker(d_I^i)$ for all $i$ (i.e., $L \cap \ker(d_I^i) \ne \{ 0 \}$ 
for any non-zero $\ZZ^n$-graded submodule $L$ of $I^i$).   
As projective resolutions, $\I$ is minimal 
if and only if so is $\J := \Phi(\I)$ 
(i.e., each $J^i$ is an essential extension of $\ker (d_J^i)$). 
Thus every $\M \in \Kom^b(\Sq)$ has a unique minimal injective resolution, 
and  any injective resolution is a direct sum of a minimal one and an exact 
complex.  For $\M \in \Kom^b(\Sq)$ and $F \subset [n]$, we define natural 
numbers $\bar{\mu}^i(F, \M)$ so that 
$$I^i \cong \bigoplus_{F \subset [n]}(S/P_F)^{\bar{\mu}^i(F,\M)},$$  
where $\I \in \Kom^b(\Sq)$ is a minimal injective resolution of $\M$. 
If $\I \in \Kom^b(\Sq)$ is a (minimal) injective resolution of 
$\M \in \Kom^b(\Sq)$,  then $\Zc(\I)$ is a (minimal) injective resolution 
of $\Zc(\M)$ in $\MMZn$.  
Hence $$\bar{\mu}^i(F, \M) = \mu^i(P_F, \Zc(\M)),$$ 
where $\mu^i(-)$ is the usual Bass number of a complex  (cf. \cite{PR}). 

Note that $\beta_i(F, -)$ and $\bar{\mu}^i(F, -)$ are 
invariants of isomorphic classes in  $D^b(\Sq)$. 

For $\M$ and $\N$, we define a complex $\Hom_S^\bullet(\M, \N)$ by 
$\Hom_S^i(\M, \N) = \prod_{j \in \ZZ} \Hom_S (M^j, N^{i+j})$ and 
the differential $d^i(f) = (\, (-1)^i f_{j+1} d_M^j + d_N^{i+j} f_j 
\,)_{j \in \ZZ}$  for $f = (f_j)_{j \in \ZZ} \in
\Hom^i_S(\M, \N)$.  Note that if $\M, \N \in \Kom^b(\MMZn)$ and each $M^i$ is
finitely generated  then $\Hom^\bullet_S(\M, \N) \in \Kom^b(\MMZn)$.

\begin{lem}\label{barmu}
Let $\I$ be a (not necessarily minimal) injective resolution 
of $\M \in \Kom^b(\Sq)$. For $F \subset [n]$, we have 
$$\bar{\mu}^i(F,\M) = \dim_k [H^i(\Hom^\bullet_S (S/P_F, \I))]_F. $$
\end{lem}

\begin{proof}
If $E^\bullet \in \Kom^b(\Sq)$ is an exact complex consisting of 
injective objects, then $E^\bullet$ splits and 
$\Hom^\bullet_S(S/P_F, E^\bullet)$ is clearly exact. 
So we may assume that $\I$ is minimal. 
Note that $\Hom_S(S/P_F, S/P_G)$ is isomorphic to 
$S/P_G$ if $F \supset G$, and 0 otherwise. Hence we have  
$[\Hom^i_S(S/P_F, \I)]_F = k^{\, \bar{\mu}^i(F,\M)}$ 
and  the differentials of 
$[\Hom_S^\bullet(S/P_F, \I)]_F$ are 0. So we are done. 
\end{proof}

\section{Functors on the Derived Category of Squarefree Modules}
Let $\Lambda$ be the incidence algebra of $2^{[n]}$. 
If $N \in \md$, then $\Hom_k(N, k)$ has a right $\Lambda$-module  
(i.e., a left $\Lambda^\op$-module) structure 
such that $(f\lambda)(a) = f(\lambda a)$ for $\lambda \in \Lambda$ and 
$a \in N$, see \cite[\S II. 3]{ARS}. 
But the opposite ring $\Lambda^\op$ of $\Lambda$ is isomorphic to 
$\Lambda$ itself by $\Lambda^\op \ni e_{F,G} \mapsto e_{G^\cmpl, F^\cmpl} 
\in \Lambda$, where $F^\cmpl := [n] \setminus F$. 
Thus $\Hom_k(-,k)$ gives a contravariant functor from $\md$ 
to itself. By the equivalence $\Sq \cong \md$ of Proposition~\ref{sq=sh}, 
we have an exact contravariant functor 
from $\Sq$ to itself. We call this functor 
the {\it Alexander duality functor}, and denote it by $\cA$. 
We have $\cA \circ \cA \cong \Id_{\Sq}$, see \cite[II, Theorem~3.3]{ARS}. 

The functor $\cA$ was defined independently by Miller~\cite{Mil} and 
R\"omer~\cite{R0} extending an idea of Eagon-Reiner~\cite{ER}. 
But their constructions of $\cA$ are different from the above one. 
R\"omer's definition is similar to ours, but it uses 
squarefree modules over an exterior algebra. 
Miller's definition uses straight modules and the Matlis duality. In fact, 
we have $\cA(M) \cong \Nc({}^*\Hom_S(\Zc(M), \E(k))(-\b1))$.   

It is easy to see that $\cA(M)_F$ is the $k$-dual of $M_{F^\cmpl}$, 
and the multiplication $\cA(M)_F \ni y \mapsto x_i y \in 
\cA(M)_{F \cup \{ i \}}$ for $i \not \in F$ is the $k$-dual of 
$M_{F^\cmpl \setminus \{ i\}} \ni y \mapsto x_i y \in M_{F^\cmpl}$. 
For example, $\cA(S(-F)) = S/P_{F^\cmpl}$ and   
$\cA(S/I_\Delta) = I_{\Delta^*}$, where 
$\Delta^* := \{ F \subset [n] \mid F^\cmpl \not 
\in \Delta \}$ is {\it (Eagon-Reiner's) Alexander dual complex} 
(\cite{ER}) of $\Delta$. 

A complex $\I \in \Kom^b(\Sq)$ is a (minimal) injective resolution of $\M$ 
if and only if the Alexander dual $\cA(\I)$ is a (minimal)
projective resolution of $\cA(\M)$. Hence we have 
$\bar{\mu}^i(F, \M) = \beta_i(F^\cmpl, \cA(\M)).$

The following is a key lemma of this section. 

\begin{lem}[{\cite[Lemma~3.20]{Y1}}]\label{Nn part}
For a squarefree module $M$ and a subset $F \subset [n]$,  
$\Nc(\Hom_S(M, \E(S/P_F)))$ is isomorphic to 
$(M_F)^* \otimes_k (S/P_F)$. Here $(M_F)^*$ is the dual $k$-vector 
space of $M_F$, but  we set the degree of $(M_F)^*$ to be 0
(since it is essentially $\Hom_k (M_F, [S/P_F]_F)$). 
In particular, $\Nc(\Hom_S(M, \E(S/P_F)))$ is squarefree. 
\end{lem}

Let $\D$ be a minimal injective resolution of $\can[n]$ in $\MMZn$ 
(according to the usual convention on dualizing complexes, we use 
$\can[n]$ instead of $\can$ itself). 
The complex $\D$ is of the form 
\begin{equation}\label{DC}
\D : 0 \too \omega^{-n} \too \omega^{-n+1} \too \cdots \too \omega^0 \too 0,
\end{equation}
$$\omega^i = \bigoplus_{\substack{F \subset [n] \\ |F| = -i}} \E(S/P_F), $$
and the differential is composed of $(-1)^{\alpha(j,F)} \cdot 
\operatorname{nat} : \E(S/P_F) \to \E(S/P_{F \setminus \{j\}})$ for $j \in F$, 
where $\operatorname{nat} : \E(S/P_F) \to \E(S/P_{F \setminus \{j\}})$ is 
induced by the natural surjection $S/P_F \to S/P_{F \setminus \{j\}}$, 
and $\alpha(j,F) := \# \{ i \in F \mid i < j \}$. 
See \cite[\S5.7]{BH}. 

\begin{prop}\label{Ext}
Let $\M \in \Kom^b(\Sq)$, and $\P \in \Kom^b(\Sq)$ its projective resolution. 
Then $\Nc(\Hom^\bullet_S (\M, \D))$, 
$\Nc(\Hom^\bullet_S(\P, \D))$ and $\Hom^\bullet_S(\P, \can[n])$ 
belong to $\Kom^b(\Sq)$, and are isomorphic in $D^b(\Sq)$.
\end{prop}

\begin{proof} By Lemma~\ref{Nn part}, 
$\Nc(\Hom^\bullet_S (\M, \D))$ and 
$\Nc(\Hom^\bullet_S(\P, \D))$ are in $\Kom^b(\Sq)$. 
Since $\Hom_S(S(-F), \can) \cong S(-F^\cmpl)$, 
$\Hom^\bullet_S(\P, \can[n])$ also belongs to $\Kom^b(\Sq)$. 
Applying \cite[Exercise~III 5.1]{GM} to $\MMZn$, 
we have $\ZZ^n$-graded quasi-isomorphisms  
$\Hom^\bullet_S(\M, \D) \to \Hom^\bullet_S(\P, \D)$ and  
$\Hom^\bullet_S(\P, \can[n]) \to \Hom^\bullet_S(\P, \D).$ 
Hence  we have quasi-isomorphisms 
$$\Nc(\Hom^\bullet_S(\M, \D)) \to \Nc(\Hom^\bullet_S(\P, \D))$$ 
and $$\Hom^\bullet_S(\P, \can[n]) = 
\Nc(\Hom^\bullet_S(\P, \can[n])) \to 
\Nc(\Hom^\bullet_S(\P, \D)).$$ 
\end{proof}

It is easy to see that 
$\cD : M^\bullet \mapsto \Nc(\Hom^\bullet_S (\M, \D))$
defines a contravariant functor from $D^b(\Sq)$ to itself. 
If $\P$ is a projective resolution of $\M$, 
$\Hom^\bullet_S(\P, \can[n])$ and $\Nc(\Hom^\bullet_S(\P, \D))$   are 
isomorphic to $\cD(\M)$ in $D^b(\Sq)$ by Proposition~\ref{Ext}. 
Hence $H^i(\cD(\M)) = \Ext_S^{n+i}(\M, \can)$ 
and $\cD \circ \cD \cong \Id_{D^b(\Sq)}$.

\begin{rem}\label{Ext_A}
Let $\Lambda$ be the incidence algebra of $2^{[n]}$. For $N \in \md$, 
the right $\Lambda$-module $\Hom_\Lambda(N,\Lambda)$ can be seen as a left 
$\Lambda$-module by the isomorphism $\Lambda^\op \cong \Lambda$ given 
in the beginning of this section. 
Similarly, $\Ext^i_\Lambda(N,\Lambda) \in \md$. 
Let $\Phi: \Sq \to \md$ be the functor of Proposition~\ref{sq=sh}, 
and let $\Pro_\Sq$ (resp. $\Pro_\Lambda$) be the full 
subcategory of $\Sq$ (resp. $\md$) consisting of projective objects. 
Then the homotopic categories $K^b(\Pro_\Sq)$ and $K^b(\Pro_\Lambda)$ 
are equivalent to $D^b(\Sq)$ and $D^b(\md)$ respectively. 
Both $\Hom_\Lambda^\bullet(\Phi(-), \Lambda)$ and 
$\Phi \circ \Hom_S^\bullet(-, \can)$ define functors from 
$K^b(\Sq) (\cong D^b(\Sq))$ to $K^b(\Pro_\Lambda) (\cong D^b(\md))$ and 
the isomorphisms 
$$\Hom_\Lambda(\Phi(S(-F)), \Lambda) \cong 
\Hom_\Lambda(\Lambda e_F, \Lambda) \cong \Lambda e_{F^\cmpl}  
\cong \Phi(\Hom_S(S(-F), \can))$$ 
give a natural equivalence 
$\Hom_\Lambda^\bullet(\Phi(-), \Lambda) \cong 
\Phi \circ \Hom_S^\bullet(-, \can)$.  
Hence if $M$ is a squarefree module, then 
$\Ext^i_\Lambda(\Phi(M),\Lambda) \cong \Phi(\Ext^i_S(M, \can))$. 
Moreover $\cD$ corresponds to the right derived functor 
$\RR \Hom_\Lambda^\bullet(-, \Lambda)$ up to translation. 

 For $N \in \md$, we have $\Hom_k(N, k) \cong 
\Hom_\Lambda(N, \bar{E})$ as left $\Lambda^\op \, (\cong \Lambda)$-modules, 
where $\bar{E}$ is the injective envelope of $\Lambda/\rad$ as a left 
$\Lambda$-module, 
see \cite[\S II. 3]{ARS}. So $\cA$ is a representable functor too. 
\end{rem}

Let $I^\bullet$ be a minimal injective resolution of 
$\M \in \Kom^b(\MMZn)$ in $\MMZn$. For $\ba \in \ZZ^n$ 
and $i \in \NN$, let $\mu^i(\m, \M)_\ba$ be the number of 
copies of $\E(S/\m)(\ba)$ which appear 
in the Krull-Schmidt decomposition of $I^i$. 

\begin{prop} 
If $\M \in D^b(\Sq)$, then $\mu^i(\m, \M)_\ba \ne 0$ 
implies $\ba$ is squarefree. Moreover 
$\mu^i(\m, \M)_F = \beta_i(F, \, \cD(\M))$
for all $F \subset [n]$. 
\end{prop}

\begin{proof}
Since $\Hom_S(S(-\ba), \E(k)) = \E(k)(\ba)$, the argument of 
\cite[Theorem~3.6]{PR} also works here.
\end{proof}

For a squarefree module $M$, we can describe 
$\cD(M) = \Nc(\Hom^\bullet(M, \D))$ explicitly. 
By Lemma~\ref{Nn part}, we have 
\begin{equation}\label{Lcomplex}
\cD(M) : 0 \too \cD^{-n}(M) \too \cD^{-n+1}(M) \too \cdots \too 
\cD^0(M) \too 0,
\end{equation}
$$\cD^i(M) = \bigoplus_{\substack{F \subset [n] \\ |F| = -i}} 
(M_F)^* \otimes_k (S/P_F).$$
As in the lemma, the degree of $(M_F)^*$ 
is $0 \in \ZZ^n$. The differential is composed of the maps 
$$(-1)^{\alpha(j,F)} \cdot (v_j)^* \otimes_k \operatorname{nat} 
: (M_F)^* \otimes_k S/P_F \to (M_{F \setminus \{j\}})^* \otimes_k 
S/P_{F \setminus \{j\}}$$ for $j \in F$. 
Here $(v_j)^*$ is the $k$-dual of the multiplication map 
$v_j: M_{F \setminus \{j\}} \ni y \mapsto x_jy \in M_F$ and ``nat" 
is the natural surjection $S/P_F \to S/P_{F \setminus \{j\}}$.  
Note that $\cD(M)$ is a complex of injective objects in $\Sq$ and it is
minimal.  Thus we have 
$$\bar{\mu}^i(F,\cD(M)) =
\begin{cases}
\dim_k M_F & \text{if $i = -|F|$,}\\
0 & \text{otherwise.}
\end{cases}$$

For a complex $\M =\{M^i, \delta^i\} \in \Kom^b (\Sq)$,  
we can also describe the complex $\cD(\M)$ in a similar way. In fact, 
$$\cD^t(\M) = \bigoplus_{i-j= t} \cD^i(M^j) 
= \bigoplus_{-|F|-j= t} (M^j_F)^* \otimes_k (S/P_F),$$ 
and the differential is given by 
$$\cD^t(\M) \supset (M^j_F)^* \otimes_k (S/P_F) \ni 
x \otimes y \mapsto d_{\cD(M^j)}(\, x \otimes y \, )+(-1)^t 
(\delta^*(x) \otimes y) \in \cD^{t+1}(\M), $$
where $\delta^* : (M^j_F)^* \to (M^{j-1}_F)^*$ is the $k$-dual of 
$\delta^{j-1}_F :M_F^{j-1} \to M_F^j$, and $d_{\cD(M^j)}$ is the 
$-|F|^{\rm th}$ differential of $\cD (M^j)$. The complex $\cD(\M)$ 
is a complex of injective objects, but it is not minimal in general.

\begin{prop}[{cf. \cite[Proposition~3.8]{Y1}}]
If $\M \in \Kom^b(\Sq)$, then 
$$\bar{\mu}^i(F,\M) = \mu^i(P_F, \Zc(\M))  
= \dim_k \, [\Ext_S^{n-|F|-i}(\M,\can)]_F.$$ 
\end{prop}

\begin{proof}
Since $\cD^2 \cong \Id_{\Sq}$, it suffices to show 
$\bar{\mu}^i(F,\cD(\M)) = \dim_k [H^{-|F|-i}(\M)]_F.$ To see this, 
we use Lemma~\ref{barmu}. The differential $d_{\cD(M^j)}$ induces 
the zero map on $[\Hom^\bullet_S(S/P_F, \cD(\M))]_F$. Thus the complex 
$[\Hom^\bullet_S(S/P_F, \cD(\M))]_F$ of $k$-vector spaces is isomorphic to 
the complex $(\M_F)^* [\, |F| \, ]$. So we are done. 
\end{proof}

The next result was proved in \cite[Theorem~2.6]{R} for the module 
case. 

\begin{cor}
If $\M \in \Kom^b(\Sq)$, then $$\beta_i(F,\M) = \dim_k
[\Ext_S^{|F|-i}(\cA(\M),\can)]_{F^\cmpl}.$$ 
\end{cor}

\begin{proof}
We have 
$\beta_i(F,\M) = \bar{\mu}^i(F^\cmpl, \cA(\M))
= \dim_k [H^{-|F^\cmpl|-i}(\cD \circ \cA(\M))]_{F^\cmpl}.$ 
\end{proof}

Let $M$ be a squarefree module. Next we will describe the complex 
$$\bP(M) := \cA \circ \cD (M)= \cA(\Nc(\Hom^\bullet(M,\D))).$$ 
For each $F \subset [n]$,  $(M_F)^\circ$ denotes a $k$-vector space with a 
bijection $\psi_F : M_F \to (M_F)^\circ$. We denote $\psi_F(y) \in 
(M_F)^\circ$ by $y^\circ$, and set $\deg(y^\circ) = 0$.  (The essential 
meaning of $M_F^\circ$ is the $k$-dual of $\Hom_k(M_F, (S/P_F)_F)$.) 
Then $$\bP^i(M) = \bigoplus_{|F| = i} (M_F)^\circ \otimes_k S(-F^\cmpl)$$ 
and the differential map is given by 
$$d(y^\circ \otimes s) = 
\sum_{j \not \in F} (-1)^{\alpha(j,F)} (x_j y)^\circ \otimes x_js.$$  

Since $\cA$ is faithful and exact, we have the following. 

\begin{cor}[cf. {\cite[Theorem~2.10]{Y}}]\label{CM <-> acyclic}
For all $i \in \ZZ$ and all $M \in \Sq$, we have 
$H^i(\bP(M)) = \cA(\Ext^{n-i}_S(M,\can))$. 
In particular, $H^i(\bP(M)) = 0$ for all $i \ne d$ 
if and only if $M$ is a Cohen-Macaulay module of dimension $d$ 
or $M=0$. 
\end{cor}

For a complex $\M =\{M^i, \delta^i \} \in \Kom^b (\Sq)$, we can 
also describe $\bP(\M) = \cA \circ \cD (\M)$ in the following way: 
$$\bP^t(\M) = \bigoplus_{i+j= t} \bP^i(M^j) 
= \bigoplus_{|F|+j= t} (M^j_F)^\circ \otimes_k S(-F^\cmpl),$$ 
and the differential is given by 
$$\bP^t(\M) \supset (M^j_F)^\circ \otimes_k S(-F^\cmpl) \ni 
y^\circ \otimes s \mapsto d_{\bP(M^j)} (y ^\circ\otimes s)+(-1)^t 
\bar{\delta}^j(y^\circ) \otimes s \in \bP^{t+1}(\M).$$ 
Here  $d_{\bP(M^j)}$ is the $|F|^{\rm th}$ differential of $\bP(M^j)$, 
and $\bar{\delta}^j : (M^j_F)^\circ \to (M^{j+1}_F)^\circ$ is induced by 
$\delta^j: M^j \to M^{j+1}$. Note that $\bP(\M)$ is 
a complex of projective objects, but not minimal in general.

Let $\P$ be a {\it minimal} projective resolution of $\M \in \Kom^b(\Sq)$.  
Thus $$P^j = \bigoplus_{F \subset [n]} S(-F)^{\beta_{-j}(F, \M)}.$$ 
For an integer $i$, the {\it $i$-linear strand} 
$P^{\bullet}_{\< i \>}$ is defined to be the complex such that 
$$P^j_{\< i \>} = 
\bigoplus_{|F|= i-j} S(-F)^{\beta_{-j}(F, \M)}$$
is a direct summand of $P^j$ and the differential $P^j_{\< i \>} \to 
P^{j+1}_{\< i \>}$ is the corresponding component of the differential 
$P^j \to P^{j+1}$ of $\P$ (so this map is  represented by a matrix of linear 
forms). The next result generalizes \cite[Theorem~4.1]{Y}. 

\begin{thm}\label{strand}
If $\M \in D^b(\Sq)$,  the $i$-linear strand $P_{\<i\>}^\bullet$ of 
$\M$ is isomorphic to  $\bP(\Ext^i_S(\cA(\M), \can))[n-i]$. 
\end{thm}

The following is immediate from Corollary~\ref{CM <-> acyclic} 
and Theorem~\ref{strand}.  

\begin{cor}[R\"omer \cite{R}]
Let $M$ be a squarefree module. Then $M$ is componentwise linear 
(i.e., the $i$-linear strand $P^\bullet_{\< i \>}$ is acyclic for 
any $i$) if and only if $\cA(M)$ is sequentially Cohen-Macaulay 
(i.e., $\Ext_S^i(\cA(M), \can)$ is a Cohen-Macaulay module of dimension 
$n-i$ for all $i$).
\end{cor}

To prove Theorem~\ref{strand}, we reconstruct  
$\P_{\<i\>}$ using the spectral sequence. 
Let $\Q$ be a (not necessarily minimal) projective resolution of 
$\M \in \Kom^b(\Sq)$. Consider the $\m$-adic filtration 
$\Q = F_0\Q \supset F_1 \Q \supset \cdots$ of 
$\Q$ with $F_i \Q = \m^i \Q$. Set $\gr (M) := \bigoplus_{i \geq 0} \m^i 
M/\m^{i+1}M$ for an $S$-module $M$, and regard it as a module over 
$\gr S = \bigoplus_{i \geq 0} \m^i/\m^{i+1} \cong S$. 
Since $Q^t$ is a free $S$-module, $Q_0^t := \bigoplus_{p+q=t}E_0^{p,q} = 
\bigoplus_{p \geq 0} \m^p Q^t / \m^{p+1} Q^t = \gr Q^t$ 
is isomorphic to $Q^t$ (if we identify $\gr S$ with $S$). 
The maps $d_0^{p,q} : E_0^{p,q} \to E_0^{p,q+1}$ make $\Q_0$ 
a cochain complex of free $\gr(S)$-modules.  
Consider the decomposition $\Q = \P \oplus C^\bullet$, where $\P$ is minimal 
and $C^\bullet$ is exact. If we identify $Q_0^t$ with 
$Q^t = P^t \oplus C^t$, the differential $d_0$ of $\Q_0$ is given by 
$(0, d_C)$. Hence we have $Q_1^t = \bigoplus_{p+q=t}E_1^{p,q} \cong P^t$. 
The maps $d_1^{p,q}: E_1^{p,q} = \m^p P^t/\m^{p+1} P^t  
\to E_1^{p+1,q} = \m^{p+1} P^{t+1}/\m^{p+2} P^{t+1}$ makes 
$\Q_1$ a cochain complex of free $\gr(S)(\cong S)$-modules whose differential 
is the ``linear term" of the differential $d_P$ of $\P$. Thus, 
under the identification $Q_1^t = P^t$, the complex 
$Q_1^\bullet$ is isomorphic to  $\bigoplus_{i \in \ZZ} \P_{\<i\>}$. 
\medskip

\noindent{\it Proof of Theorem~\ref{strand}.} 
Since  $\cA \circ \cD \circ \cD \circ \cA \cong 
\Id_{D^b(\Sq)}$, it suffices to prove the $i$-linear strand of 
$\cA \circ \cD (\M)$ is isomorphic to $\bP(H^{-n+i}(\M))[n-i]$. 
Recall that $\bP(\M) =\cA \circ \cD(\M)$ is a complex of projective objects. 
Set $\Q = \bP(\M)$, and consider the $\m$-adic filtration $F_i \Q = \m^i\Q$ 
of $\Q$. Under the above notation, the differential  
$d_0^t : Q_0^t \cong \bP^t(\M) \to Q_0^{t+1} \cong 
\bP^{t+1}(\M)$ is given by $(-1)^t \delta$. Thus  
$$Q_1^t \cong \bigoplus_{|F|+j=t} H^j(\M) \otimes_k S(-F^\cmpl)
= \bigoplus_{l+j=t}\bP^l(H^j(\M)),$$ 
and the differential of $\Q_1$ is induced by that of 
$\bP(M^j)$. Hence we can easily check that $\Q_1$ is isomorphic to 
$\bigoplus_{j \in \ZZ} \bP(H^j(\M))[-j]$. 
By the remark before this proof, the $i$-linear strand of $\cA \circ \cD 
(\M)$ is isomorphic to $\bP(H^{-n+i}(\M))[n-i]$. 
\qed

\begin{thm}\label{DADADA} 
We have a natural equivalence 
$\cD \circ \cA \circ \cD \circ \cA \circ \cD \circ \cA \cong 
\bT^{2n}$ in $D^b(\Sq)$, where $\bT$ is the translation functor 
(i.e., $\bT^{2n} : \M \mapsto \M[2n]$). 
\end{thm}

\begin{proof}
For $\M = \{M^i, \delta^i \} \in \Kom^b(\Sq)$, the complex 
$\Hom^\bullet_S(\bP(\M),\can[n])$ is isomorphic to 
$\cD\circ \cA \circ \cD \, (\M)$ in $D^b(\Sq)$. We have  
\begin{eqnarray*}\Hom^i_S(\bP(\M),\can[n]) 
&=& \Hom_S \, (\bigoplus_{-i-n = |F|+j } 
(M^j_F)^\circ \otimes_k S(-F^\cmpl), \ \can \ ) \\
&=& \bigoplus_{-i-n = |F|+j } 
(M^j_F)^* \otimes_k S(-F)\\
&=& \bigoplus_{i = -n-|F|+j} (M^{-j}_F)^* \otimes_k S(-F). 
\end{eqnarray*} 
Here we simply denote the dual vector space of $(M^{-j}_F)^\circ$ by 
$(M^{-j}_F)^*$, since $(M^{-j}_F)^\circ \cong M^{-j}_F$ as $k$-vector 
spaces (only the degrees are different). 
Also here  $\deg (M^{-j}_F)^* =0 \in \ZZ^n$. 
The differential of $\Hom^\bullet_S(\bP(\M),\can[n])$ is 
given by 
$$(M^{-j}_F)^* \otimes_k S(-F) \ni y \otimes s \mapsto 
\sum_{l \in F} (-1)^{\alpha(l,F)+n+|F|-j} \, 
v_l^*(y) \otimes x_ls + (-1)^{n-1}\delta^*(y)\otimes s,$$
where $v_l^* : (M_F^{-j})^* \to 
(M_{F \setminus \{l \}}^{-j})^*$ is the $k$-dual of  
$v_l: M_{F\setminus \{ l \}}^{-j} \ni z \mapsto x_l z \in M_F^{-j}$, 
and $\delta^* : (M^{-j}_F)^* \to (M^{-j-1}_F)^*$ is the $k$-dual of 
$\delta^{-j-1} : M^{-j-1}_F \to M^{-j}_F$. 

Similarly, $\bP(\cA(\M))$ represents 
$\cA \circ \cD \circ \cA \, (\M)$ in $D^b(\Sq)$, and we have 
\begin{eqnarray*}
\bP^i(\cA(\M))  
&=& \bigoplus_{i=|F|+j} 
(\cA(M^j)_F)^\circ \otimes_k S(-F^\cmpl)\\
&=& \bigoplus_{i=|F|+j} 
(M^{-j}_{F^\cmpl})^* \otimes_k S(-F^\cmpl)\\
&=& \bigoplus_{i=n-|F|+j} 
(M^{-j}_F)^* \otimes_k S(-F). 
\end{eqnarray*}
Also here, we simply denote  
$(\cA(M^{-j})_F)^\circ= ((M^{-j}_{F^\cmpl})^*)^\circ$ by
$(M^{-j}_{F^\cmpl})^*$. The differential of the above complex is
given by  $$(M^{-j}_F)^* \otimes_k S(-F) \ni y \otimes s \mapsto 
\sum_{l \in F} (-1)^{\alpha(l,F^\cmpl)} \, 
v_l^*(y) \otimes x_ls + (-1)^{|F^\cmpl|+j} \, \delta^*(y)\otimes s.$$

For an integer $l \in \ZZ$, set $\beta(l) := 1$ if 
$l \equiv 1,2 \pmod{4}$, and $\beta(l) := 0$ if $l \equiv 3,0 \pmod{4}$. 
We also set $\alpha(A,B) := \# \{(a,b) \, | \, a > b, a \in A, b \in B \}$ 
for $A, B \subset [n]$.
Then the multiplication by $(-1)^{\alpha(F, [n])+\beta(|F|-j)+|F|n + j}$
on $(M^{-j}_F)^* \otimes_k S(-F)$, which can be regarded as a submodule of 
both $\Hom^{-n-|F|+j}_S(\bP(\M),\can[n])$ and $\bP^{n-|F|+j}(\cA(\M))$, 
induces quasi-isomorphism between $\Hom^\bullet_S(\bP(\M),\can[n])$ and
$\bT^{2n} \circ \bP(\cA(\M))$. So $\cD \circ \cA \circ \cD \cong 
\bT^{2n} \circ \cA \circ \cD \circ \cA$ as a functor on $D^b(\Sq)$. Since 
$(\cA \circ \cD \circ \cA) \circ (\cA \circ \cD \circ \cA) \cong 
\Id_{D^b(\Sq)}$, we get the assertion. 
\end{proof}

\begin{exmp}  
For $F \subset [n]$, we have the following. 
\begin{eqnarray*}
& & \cD \circ \cA \circ \cD \circ \cA \circ \cD \circ \cA \, (S(-F))\\
&&= \cD \circ \cA \circ \cD \circ \cA \circ \cD \, (S/P_{F^\cmpl})\\
&&= \cD \circ \cA \circ \cD \circ \cA \, 
((S/P_{F^\cmpl}) (-F^\cmpl)[-|F|+n])\\
&&= \cD \circ \cA \circ \cD \, ((S/P_F)(-F)[|F|-n])\\
&&= \cD \circ \cA \, ((S/P_F)[n])\\
&&= \cD \, (S(-F^\cmpl)[-n])\\
&&= S(-F)[2n]. 
\end{eqnarray*}
\end{exmp}

\section{Relation to Koszul duality} 
Let $S = k[x_1, \ldots, x_n]$ be a polynomial ring as in 
the previous sections, and $E :=\bigwedge S_1^* = 
k \<e_1, \ldots, e_n \>$ an exterior algebra. 
$E$ is a $\ZZ^n$-graded ring with 
$\deg (e_i) = (0, \ldots, 0, -1, 0, \ldots, 0) = - \deg (x_i)$ 
where $-1$ is in the $i^{\rm th}$ position. When we regard $S$ and $E$ 
as  $\ZZ$-graded rings, we set $\deg(x_i)=1$ and $\deg (e_i) = -1$ for all 
$i$. In this paper, $E$-modules are {\it left} $E$-modules unless otherwise 
specified. For a $\ZZ^n$-graded $E$-module $M$ and $\ba \in \ZZ^n$, 
$M_\ba$ means the degree $\ba$ component of $M$, and $M(\ba)$ is the shifted 
module with $M(\ba)_\bb = M_{\ba+\bb}$ as in the polynomial ring case. 

Denote the category of finitely generated $\ZZ$-graded $S$-modules 
(resp. $E$-modules) by $\mmS$ (resp. $\mmE$). Although $\mmS$ and $\mmE$ 
are far from equivalent, a famous theorem of 
Bernstein-Gel'fand-Gel'fand~\cite{BGG} 
states that $D^b(\mmS) \cong D^b(\mmE)$ as triangulated categories. 
First, we will see that this 
equivalence also holds in the $\ZZ^n$-graded context. 
Denote the category of finitely generated $\ZZ^n$-graded $S$-modules 
(resp. $E$-modules) by $\mmSn$ (resp. $\mmEn$).

There are several papers concerning the Bernstein-Gel'fand-Gel'fand 
correspondence. But their conventions are not quite the same. 
In this paper, we basically follow \cite{ES}, 
which is well suited for our purpose. 
Here we give functors defining $D^b(\mmSn) \cong D^b(\mmEn)$. 
For $M \in \mmSn$, we define  $\cR (M)= \Hom_k(E(-{\bf 1}),M)$ to be  
a $\ZZ^n$-graded cochain complex of free $E$-modules as follows. 
(The original definition is $\cR (M)= \Hom_k(E,M)$, 
but we use this grading. We will also shift the grading of $\cL(N)$ 
defined below.)  We can define a $\ZZ^n$-graded left $E$-module 
structure on $\Hom_k(E(-{\bf 1}),M_\ba)$ by $(a f)(e) = f(ea)$. 
Then $\Hom_k(E(-{\bf 1}),M_\ba) \cong E(-\ba)^{\oplus \dim_k M_\ba}$.  
Set the cohomological degree of $\Hom_k(E(-{\bf 1}), M_\ba)$ to be 
$||\ba||:= \sum_{j \in [n]} a_j$.  The differential of $\cR(M)$ is defined by 
$$\Hom_k(E(-{\bf 1}),M_\ba) \ni f \mapsto [e \mapsto \sum_{i \in [n]} 
x_i f(e_ie)] \in \bigoplus_{i \in [n]} \Hom_k(E(-{\bf 1}), 
M_{\ba + \varepsilon_i}),$$
where $\varepsilon_i \in \NN^n$ is the squarefree vector whose support is 
$\{ i \}$. We also define the complex $\cR(\M) = \bigoplus_{j \in  \ZZ} 
\Hom_k(E(-{\bf 1}), M^j)$ for a complex $\M = \{ M^j, \delta^j\}$ in $\mmSn$. 
The cohomological degree $i$ component of $\cR(\M)$ is 
$\bigoplus_{i = j+ ||\ba||} \Hom_k(E(-{\bf 1}), M^j_\ba)$ 
and  the differential is given by 
$$\cR^i(\M) \supset \Hom_k(E(-\b1), M^j_\ba) \ni f \mapsto 
 d_{\cR(M^j)}(f) + (-1)^i (\delta^j \circ f) \in 
\cR^{i+1}(\M),$$ where $d_{\cR(M^j)}$ is 
the $||\ba||^{\rm th}$ differential of $\cR(M^j)$.  
We can apply $\cR$ to a $\ZZ$-graded complex $\M \in \Kom^b(\mmS)$
(in this case, we replace $S(-\b1)$ by $S(-n)$).  
Then $\cR$ is equivalent to the functor given in 
\cite{BGG, BGS, ES} up to degree shifting. For $\M \in \Kom^b(\mmSn)$, 
$\cR(\M)$ has only finitely many non-vanishing cohomologies.  
And $\cR$ induces a covariant functor from $D^b(\mmSn)$ to $D^b(\mmEn)$, 
which is also denoted by $\cR$. 

Next, we will define the functor $\cL : \Kom^b(\mmEn) \to \Kom^b(\mmSn)$. 
Set $\cL(\N) = \bigoplus_{i \in \ZZ} S(-{\bf 1}) \otimes_k N^i$ 
for a complex $\N = \{N^i, \delta^i\}$ in $\mmEn$. The cohomological 
degree of $\cL(\N)$ is given by $\cL^i(\N)= \bigoplus_{i=j-||\ba||} 
S(-{\bf 1}) \otimes_k N_\ba^j$.  And the differential is defined by 
$$\cL^i(\N) \supset S(-{\bf 1}) \otimes_k N_{\ba}^j \ni s \otimes y \mapsto 
\sum_{l \in [n]} x_l s \otimes e_l y + (-1)^i (s \otimes \delta^j(y)) 
\in \cL^{i+1}(\N).$$ 
If we apply $\cL$ to $\ZZ$-graded complexes, 
it is equivalent to the functor given in \cite{BGG, BGS, ES} up to degree 
shifting. If $\N$ is bounded, so is $\cL(\N)$. And $\cL$
induces a covariant functor from $D^b(\mmEn)$ to $D^b(\mmSn)$, 
which is also denoted by $\cL$.  

In the $\ZZ$-graded case, 
Bernstein-Gel'fand-Gel'fand~\cite{BGG} (see also \cite{BGS,ES}) states that  
$\cL : \Kom^b(\mmE) \to \Kom^b(\mmS)$ is a left adjoint 
to $\cR : \Kom^b(\mmS) \to \Kom^b(\mmE)$, that is, 
we have a natural isomorphism 
$$\varphi: \Hom_{\Kom^b(\mmS)}(\cL(\N), \M) \,  \stackrel{\cong}
{\too} \, \Hom_{\Kom^b(\mmE)}(\N, \cR(\M))$$  
for $\M \in \Kom^b(\mmS)$ and $\N \in \Kom^b(\mmE)$.  
Moreover, the map $\cL \circ \cR (\M) \to \M$ associated to 
the identity map $\cR(\M) \to \cR(\M)$ is a quasi-isomorphism.   
Similarly, the map $\N \to \cR \circ \cL (\N)$ associated to 
the identity map $\cL(\N) \to \cL(\N)$ is a quasi-isomorphism. Hence 
$\cR$ and $\cL$ define an equivalence $D^b(\mmS) \cong D^b(\mmE)$. 

We can regard a $\ZZ^n$-graded module as a $\ZZ$-graded module 
by $M_i = \bigoplus_{||\ba|| = i} M_\ba$. In this sense,  
$\Kom^b(\mmSn)$ and $\Kom^b(\mmEn)$ are (non-full) subcategories of 
$\Kom^b(\mmS)$ and $\Kom^b(\mmE)$ respectively. 
If $\M \in \Kom^b(\mmSn)$ and $\N \in \Kom^b(\mmEn)$, then the restriction 
of $\varphi$ gives the isomorphism 
$$\Hom_{\Kom^b(\mmSn)}(\cL(\N), \M) \,  \cong 
\, \Hom_{\Kom^b(\mmEn)}(\N, \cR(\M)).$$  
Thus the quasi-isomorphisms $\cL \circ \cR (\M) \to \M$ and 
$\N \to \cR \circ \cL (\N)$ are $\ZZ^n$-graded. 
Hence we have the following. 

\begin{thm}[BGG correspondence ($\ZZ^n$-graded version)]\label{BGG} 
The functors $\cR$ and $\cL$ define 
an equivalence of triangulated categories $D^b(\mmSn) \cong D^b(\mmEn)$. 
\end{thm}

The functors $\cR$ and $\cL$ are closely related to $\cD$ and $\cA$  
of the previous section. To see this, we recall the definition of 
a squarefree module over $E$. 

\begin{dfn}[R\"omer~\cite{R0}]
A $\ZZ^n$-graded $E$-module $N = \bigoplus_{\ba \in \ZZ^n} N_\ba$ is 
{\it squarefree} if $N$ is finitely generated and 
$N = \bigoplus_{F \subset [n]} N_{-F}$. 
\end{dfn}

For example, a monomial ideal of $E$ is always squarefree. 
We denote the full subcategory of $\mmEn$ consisting of all the squarefree 
$E$-modules by $\SqE$. We have the functors $\cS: \SqE \to \SqS$ 
and $\cE: \SqS \to \SqE$ giving an equivalence $\SqS \cong \SqE$. Here 
$\cS(N)_F = N_{-F}$ for $N \in \SqE$, and the multiplication map 
$\cS(N)_F \ni y \mapsto x_i y \in \cS(N)_{F \cup \{ i \}}$ for $i \not \in F$ 
is given by  $\cS(N)_F =N_{-F} \ni z \mapsto (-1)^{\alpha(i, F)} e_i z 
\in N_{-(F \cup \{ i \})}= \cS(N)_{F \cup \{ i \}}$. 
See \cite{R0} for further information.  

We have $D^b(\Sq_S) \cong D^b_{\Sq_S}(\mmSn)$ 
by \cite[Exercises~III 2.2 ]{GM} (but use projective resolutions 
instead of injective resolutions).
So $D^b(\Sq_S)$ can be seen as a full subcategory of $D^b(\mmSn)$. 
On the other hand, for $N \in \mmEn$, set 
$$N' := \bigoplus_{\ba \in \NN^n} N_{-\ba} \quad \text{and} \quad 
N'' := \bigoplus_{\substack{\ba \in \NN^n \text{and $\ba$ is } \\
\text{{\it not} squarefree}}} 
N_{-\ba} \subset N'. $$
Note that $N'$ and $N''$ are $E$-modules, and   
$\Nc(N) := N'/N''$ is squarefree. If all cohomologies of  
$\N \in \Kom^b (\mmEn)$ are squarefree, then 
$\Nc(\N)$ and $\N$ are isomorphic in $D^b(\mmEn)$. 
Hence we have $D^b_\SqE(\mmEn) \cong D^b(\SqE)$.

Comparing $\cL$ and $\bP = \cA \circ \cD$ defined 
in the last section, we have the following.  

\begin{prop}\label{LE} 
If $\N$ is a (bounded) complex of squarefree $E$-modules, then 
$\cL(\N) = S(-{\bf 1}) \otimes_k \N$ is a (bounded) complex of squarefree 
$S$-modules. Hence $\cL$ gives a functor from $D^b(\SqE)$ to 
$D^b(\SqS)$. Moreover, for $\M \in \Kom^b(\SqS)$, we have 
$\cL \circ \cE (\M) = \cA \circ \cD (\M)$. 
\end{prop}

On the other hand, $\cR(M)$ is not a complex of squarefree $E$-modules. 
In fact, a free $E$-module $E(-\ba)$ is not squarefree unless $\ba =0$. 
But we have the following.

\begin{prop}\label{SR}
If $M^\bullet \in D^b(\SqS)$, 
then $\cR(\M) \in D^b_\SqE(\mmEn) \cong D^b(\SqE)$. 
Moreover, we have a natural equivalence $\cS \circ \cR \cong \cD \circ \cA$.  
\end{prop}

\begin{proof} We have $\M \cong \cA \circ \cD \circ \cD \circ \cA (\M) = 
\cL \circ \cE \circ \cD \circ \cA (\M)$ in 
$D^b(\SqS)$ (and in $D^b(\mmSn)$) by Proposition~\ref{LE}. 
From Theorem~\ref{BGG}, 
$\cR(\M) \cong \cR \circ \cL \circ \cE \circ \cD \circ 
\cA (\M) \cong \cE \circ \cD \circ \cA (\M) \in D^b_{\SqE}(\mmEn)$. 
Since $\cS \circ \cE \cong \Id_{D^b(\SqS)}$, we are done. 
\end{proof}

Let $R = \bigoplus_{i \geq 0} R_i$ be an $\NN$-graded associative $k$-algebra 
such that $\dim_k R_i < \infty$ for all $i$ and $R_0 \cong k^m$ for some 
$m \in \NN$ as an algebra. Then $\rad := \bigoplus_{i > 0} R_i$ is the graded 
Jacobson radical.  We say $R$ is {\it Koszul}, 
if a left $R$-module $R/\rad$ admits a graded projective resolution 
$$\cdots \to P^{-2} \to P^{-1} \to P^0 \to R/\rad \to 0$$ such that $P^{-i}$ 
is generated by its degree $i$ component, that is, $P^{-i} = R P_i^{-i}$ 
(we say such a resolution is a {\it linear} resolution). 
If $R$ is Koszul, it is a quadratic ring, and its {\it quadratic dual 
ring} $R^!$ (see \cite[Definition~2.8.1]{BGS}) 
is Koszul again, and isomorphic to the opposite ring of the Yoneda algebra 
$E(R) := \bigoplus_{i \geq 0}\Ext^i_R(R/\rad,R/\rad)$. 

Let $\mdR$  be the category 
of finitely generated $\ZZ$-graded left $R$-modules. 
If $R$ is a Koszul algebra with $R_i = 0$ for $i \gg 0$, 
and $R^!$ is left noetherian, we have functors  
$$DF : D^b(\mdR) \ni \N \mapsto R^! \otimes_{R_0} \N \in D^b(\mmR)$$ and 
$$DG : D^b(\mmR) \ni \M \mapsto  \Hom_{R_0} (R, \M) \in D^b(\mdR)$$   
giving the equivalence $D^b(\mdR) \cong D^b(\mmR)$ 
called {\it Koszul duality}, see \cite[Theorem~2.12.6]{BGS}. 
The exterior algebra $E$ is a Koszul algebra with $E^! \cong S$. 
Thus the Bernstein-Gel'fand-Gel'fand correspondence is a 
classical example of Koszul duality. 

Let $\Lambda$ be the incidence algebra of $2^{[n]}$ over $k$. Then $\Lambda$ 
has an $\NN$-grading with $\deg(e_{F,G}) = |F \setminus G|$. 
Note that $\Lambda_0 = \bigoplus_{F \subset [n]} k e_F   \cong k^{2^n}$. 
 For each $i \in \ZZ$, let $\mdZ(i)$ be the full subcategory of $\mdZ$ 
consisting of $N \in \mdZ$ such that $N_j = \bigoplus_{|F| = i+j} N_F$ 
for all $j \in \ZZ$. The forgetful functor gives an equivalence 
$\mdZ(i) \cong \md$ for all $i \in \ZZ$, and $D^b(\mdZ(i)) \cong 
D^b_{\mdZ(i)}(\mdZ)$ is a full subcategory of $D^b(\mdZ)$.   
Since $\Lambda / \rad \Lambda \cong \bigoplus_{F \subset [n]} k e_F$ as left 
$\Lambda$-modules, and each $k e_F$ has a linear projective resolution 
$$\cdots \to \bigoplus_{\substack{G \supset F \\|G| = |F|+2}} \Lambda e_G 
\to \bigoplus_{\substack{G \supset F \\|G| = |F|+1}} \Lambda e_G 
\to \Lambda e_F \to k e_F \to 0$$ 
(here we regard $k e_F$ and $\Lambda e_G $ as objects in $\mdZ(|F|)$), 
$\Lambda$ is Koszul. 
To see this we can use Proposition~\ref{sq=sh}. 
In fact, a minimal free resolution of a squarefree module $(S/P_F)(-F)$, 
which corresponds to $k e_F$,  is given by the Koszul complex 
with respect to $\{x_i \mid i \not \in F \}$. 

\begin{lem}\label{self dual}
The quadratic dual ring $\Lambda^!$ of $\Lambda$ is 
isomorphic to $\Lambda$ itself. 
\end{lem}

One might think $\Lambda^!$ should be a ``negatively graded ring", 
since $\Lambda^!$ is generated by $\Hom_{\Lambda_0}(\Lambda_1, \Lambda_0)$ 
as a $\Lambda^!_0$-algebra. But we use the same convention as \cite{BGS} 
here, so we regard $\Lambda^!$ as a positively graded ring with 
$\Lambda^!_1 = \Hom_{\Lambda_0}(\Lambda_1, \Lambda_0)$.

\begin{proof}
Let $T := T_{\Lambda_0} \Lambda_1 = \Lambda_0 \oplus \Lambda_1 \oplus 
(\Lambda_1 \otimes_{\Lambda_0} \Lambda_1) \oplus \cdots 
= \bigoplus_{i \geq 0} \Lambda_1^
{\otimes i}$ be the tensor ring of $\Lambda_1 = 
\< \, e_{F \cup \{ i \}, F} \mid  F \subset [n], \, i \not \in F \, \>$. 
(See \cite[\S~2.7]{BGS} 
for the linear algebra over a semisimple algebra $\Lambda_0$ used here.)
Then $\Lambda \cong T/I$, where 
$$I = (\,  e_{F \cup \{ i, j\}, F \cup \{ i \}} \otimes 
e_{F \cup \{ i \}, F}  - e_{F \cup \{ i, j\}, F \cup \{ j \}}
\otimes   e_{F \cup \{ j \}, F} 
\mid F \subset [n], i, j \not \in F \, )$$ is a two sided ideal.  
Let $\Lambda_1^* := \Hom_{\Lambda_0}(\Lambda_1, \Lambda_0)$ 
be the dual of the {\it left} $\Lambda_0$-module $\Lambda_1$. 
Then  $\Lambda_1^*$ has a right $\Lambda_0$-module 
structure such that $(fa)(v) = (f(v))a$, and a left $\Lambda_0$-module 
structure such that $(af)(v) =  f(va)$, where $a \in \Lambda_0$, 
$f \in \Lambda_1^*$, $v \in \Lambda_1$. As a left (or right) 
$\Lambda_0$-module, $\Lambda_1^*$ is generated by  
$\{ \, e_{F, F \cup \{ i \}}^* \mid  F \subset [n], 
\, i \not \in F \,\}$, where 
$$e_{F, F \cup \{ i \}}^*(e_{G \cup \{ j \},G}) = 
\delta_{F,G} \, \delta_{i,j} \, e_{F \cup \{ i \}}.$$ 
Let $T^* = T_{\Lambda_0} \Lambda_1^*$ be the tensor ring of $\Lambda_1^*$. 
Note that $e^*_{F, F \cup \{ i \}} \otimes e^*_{G, G \cup \{ j \}} 
\ne 0$ if and only if $F \cup \{ i \} = G$. 
We have that $(\Lambda_1^* \otimes_{\Lambda_0} \Lambda_1^*)$ 
is isomorphic to $(\Lambda_1 \otimes_{\Lambda_0} \Lambda_1)^* = 
\Hom_{\Lambda_0}(\Lambda_1 \otimes_{\Lambda_0} \Lambda_1, \Lambda_0)$ 
via $(f \otimes g)(v \otimes w) = g(vf(w))$, where 
$f, g \in \Lambda_1^*$ and $v, w \in \Lambda_1$. 
In particular, 
$$( e^*_{F, F \cup \{ i \}} \otimes e^*_{F \cup \{ i \}, F \cup \{ i, j \}}) 
( e_{F \cup \{ i, j \}, F \cup \{ i \}} \otimes e_{F \cup \{ i \}, F} )
= e_{F \cup \{ i, j\}}.$$
Easy computation shows that the quadratic dual ideal 
$$I^\bot = (\, f \in \Lambda_1^* \otimes \Lambda_1^* \mid \text{$f(v) = 0$ 
for all $v \in I_2 \subset  \Lambda_1 \otimes \Lambda_1 = T_2$} \,) 
\subset T^*$$ of $I$ is equal to  
$$(\,  
e^*_{F, F \cup \{ i \}} \otimes e^*_{F \cup \{ i \}, F \cup \{ i,j \}} + 
 e^*_{F, F \cup \{ j \}}  \otimes e^*_{ F \cup \{ j \}, F \cup \{ i, j \}}
 \mid F \subset [n], i, j \not \in F, i \ne j \, ).$$ 
The  $k$-algebra homomorphism defined by 
$$\Lambda_0 \ni e_F \mapsto e_{F^\cmpl} \in \Lambda_0^! \, (= \Lambda_0)
\quad \text{and} \quad \Lambda_1 \ni e_{F \cup \{ i\},F} 
\mapsto (-1)^{\alpha(i, F)} 
e^*_{(F \cup \{ i\})^\cmpl, F^\cmpl} \in \Lambda^!_1$$
gives a graded isomorphism $\Lambda \cong \Lambda^!$. 
\end{proof}

Since $\Lambda \, (\cong \Lambda^!)$ is an artinian algebra, 
we have the functors $DF$ and $DG$ 
defining $D^b(\mdZ) \cong D^b({\bf gr.mod}_{\Lambda^!})$. 
In the next result, we will denote the contravariant functors 
from $D^b(\md)$ to itself induced by $\cD$ and $\cA$ for $D^b(\Sq)$ 
(under the equivalence $\Sq \cong \md$ of Proposition~\ref{sq=sh}) 
also by $\cD$ and $\cA$. 

\begin{thm}\label{self Koszul dual}
Let the notation be as above. If $\N \in D^b(\mdZ(0))$,  then 
we have $DF(\N) \in D^b(\mdZ(n))$ and $DG(\N) \in D^b(\mdZ(-n))$  
under the isomorphism $\Lambda^! \cong \Lambda$ of Lemma~\ref{self dual}.
By the equivalence $\mdZ(j) \cong \md$, 
$DF$ and  $DG$ give endofunctors of $D^b(\md)$. Then  
$DF \cong \cA \circ \cD$ and $DG \cong \cD \circ \cA$ as 
endofunctors of $D^b(\md)$. 
\end{thm}

\begin{proof}
First, we recall the construction of 
$DF : D^b(\mdZ) \ni \N \mapsto \Lambda^! \otimes_{\Lambda_0} \N \in 
D^b({\bf gr.mod}_{\Lambda^!})$ 
under the same notation as the proof of the previous lemma. 
Note that  $\Lambda^!_0 = \Lambda_0 = \bigoplus_{ F \subset [n]} k e_F$. 
The component $(DF)^t(\N)$ of cohomological degree $t$  
is $\bigoplus_{t=i+j} \Lambda^! \otimes_{\Lambda_0} 
N_j^i$.  For $N \in \mdZ$, a left $\Lambda^!$-module 
$\Lambda^! \otimes_{\Lambda_0} N = 
\bigoplus_{F \subset [n]} \Lambda^! e_F \otimes_k N_F$ 
is generated by $\{\, e_F \otimes n_F \mid \text{$F \subset [n]$ and 
$n_F \in N_F$} \, \}$.  
If $N \in \mdZ(0)$, the degree of $e_F \otimes n_F$ is 
$\deg(e_F) - \deg(n_F) = -|F|$.  For $\N = \{N^i, \delta^i\} 
\in \Kom^b(\mdZ)$, the differential of $DF(\N)$ is given by 
$$(DF)^t(\N) \ni e_F \otimes n_F \mapsto 
(-1)^t \sum_{l \not \in F} e^*_{F, F \cup \{ l \}} \otimes 
(\, e_{F \cup \{ l \}, F} \cdot n_F \, ) + e_F \otimes \delta(n_F),$$
see \cite[Theorem~2.12.1]{BGS}. 

The graded isomorphism $\Lambda \stackrel{\cong}{\to} \Lambda^!$ 
makes $M \in {\bf gr.mod}_{\Lambda^!}$ 
a graded left $\Lambda$-module (without changing the grading of $M$), 
and gives an equivalence ${\bf gr.mod}_{\Lambda^!} \cong 
{\bf gr.mod}_\Lambda$. From now on, we regard $DF$ as an endofunctor of 
$D^b(\mdZ)$ by the equivalence ${\bf gr.mod}_{\Lambda^!} \cong 
{\bf gr.mod}_\Lambda$. 
So we have $DF(N) = \bigoplus_{F \subset [n]} \Lambda e_{F^\cmpl} 
\otimes_k N_F$ for $N \in \mdZ$. If $N \in \mdZ(0)$, then the degree of 
$e_{F^\cmpl} \otimes n_F \in \Lambda e_{F^\cmpl} \otimes_k N_F \subset DF(N)$ 
is $-|F|= |F^\cmpl| - n$. Thus $DF(N) \in \mdZ(n)$. For $\N \in D^b(\mdZ(0))$, 
the cohomological degree of $DF(\N)$ is given by 
 $(DF)^t(\N)= \bigoplus_{t=j+|F|} \Lambda e_{F^\cmpl} \otimes_k N^j_F$, 
and the differential sends $e_{F^\cmpl} \otimes n_F \in (DF)^t(\N)$ to  
$$  \sum_{l \not \in F} (-1)^{t + \alpha(l,F)} \,  
e_{F^\cmpl, (F \cup \{ l \})^\cmpl} \otimes 
(\, e_{F \cup \{ l \}, F} \cdot n_F \, ) + e_{F^\cmpl} \otimes 
\delta(n_F).$$

In Section~3, we study the endofunctor $\bP = \cA \circ \cD$ on $D^b(\Sq)$. 
Under the equivalence $\Sq \cong \md$ of 
Proposition~\ref{sq=sh}, this functor induces an endofunctor of $D^b(\md)$. 
We also denote it by $\cA \circ \cD$. Then for $\N \in D^b(\md)$, 
the component $(\cA \circ \cD)^t (\N)$ of cohomological degree $t$ is 
$\bigoplus_{t = j+|F|} \Lambda e_{F^\cmpl} \otimes_k N^j_F$, and 
an element $e_{F^\cmpl} \otimes n_F \in \Lambda e_{F^\cmpl} \otimes_k N_F 
\subset(\cA \circ \cD)^t (\N)$  is  sent to  
$$ \sum_{l \not \in F} (-1)^{\alpha(l,F)} \,  
e_{F^\cmpl, (F \cup \{ l \})^\cmpl} \otimes 
(\, e_{F \cup \{ l\}, F}  \cdot  n_F \, ) + 
(-1)^t \, e_{F^\cmpl} \otimes \delta(n_F)$$ by the differential. 
A quasi-isomorphism  
$(DF)^t(\N) \ni x \mapsto (-1)^{\beta (t-1)} x \in (\cA \circ \cD)^t (\N)$ 
gives a natural equivalence $DF \cong \cA \circ \cD$, where $\beta(-)$ is 
the function defined in the proof of Theorem~\ref{DADADA}. The natural 
equivalence $DG \cong \cD \circ \cA$ can be proved in a similar way. 
\end{proof}

\section{Local Cohomology modules as holonomic $D$-modules}
In this section, we study a local cohomology module $H_{I_\Delta}^i(S)$. 
The following result was essentially obtained by 
\Mus~\cite{Mus} and Terai~\cite{T}, and can be proved by the 
same argument as the proof of \cite[Theorem~2.11]{Y1}. 

\begin{thm}[cf. \cite{Mus,T,Y1}]
Let $\Gamma_{I_\Delta}$ be the local cohomology functor with supports in 
$I_\Delta$. Then $\Gamma_{I_\Delta}(\D) \in D^b(\Str)$ and 
$\Zc \circ \cD(S/I_\Delta) \cong \Gamma_{I_\Delta}(\D)$.  
In particular, $H_{I_\Delta}^i(S)(-\b1) \cong H_{I_\Delta}^i(\can) 
\cong \Zc(\Ext_S^i(S/I_\Delta, \can))$. 
\end{thm}

See \cite{Mil2, Y3} for further results on minimal flat resolutions 
of $\Gamma_{I_\Delta}(\D)$. 

\medskip

In the rest of this section, we assume that $\chara(k) =0$. 
Let $$A := A_n(k) = k\<x_1, \ldots, x_n, \partial_1, \ldots, \partial_n \>$$ 
be the Weyl algebra acting on $S$, and let 
$\{ F_i \}_{ i \geq 0}$ 
with $F_i = \< \, x^\ba \partial^\bb \mid |\ba| + |\bb| \leq i \, \>$ 
be the Bernstein filtration of $A$. 
Here $|\ba| = \sum_{i=1}^n |a_i|$ for $\ba = (a_1, \ldots, a_n)$. 
Then the associated graded ring  
$\grF A := \bigoplus_{i \geq 0} F_i/F_{i-1}$ is isomorphic to the polynomial 
ring $k[\bar{x}_1, \ldots, \bar{x}_n, \bar{\partial}_1, \ldots, 
\bar{\partial}_n]$ of $2n$ variables. See, for example, \cite{Bj}. 

In \cite{Y1}, the author pointed out that a straight $S$-module $M$ has a 
holonomic $A$-module structure. But if we consider the $\ZZ^n$-grading, 
the {\it left} $A$-module structure given in \cite{Y1} is 
somewhat unnatural. So we will give a more natural treatment here. 

Let $M$ be a left $A$-module. 
Set 
\begin{equation}\label{left}
M_\rat : = \bigoplus_{\ba \in \ZZ^n} M_\ba, 
\quad 
\text{where} \quad  M_\ba = \{ \, y \in M \mid 
\text{$(x_i\partial_i) \, y =  a_iy$ for all $i$} \, \}. 
\end{equation}
Then $M_\rat$ is an $A$-submodule  
with $x_i M_\ba \subset M_{\ba +\varepsilon_i}$ and  
$\partial_i M_\ba \subset M_{\ba-\varepsilon_i}$. 
In particular, $M_\rat$ is a $\ZZ^n$-graded $S$-module. 
For example, $S_\rat = S$ and the $\ZZ^n$-grading given 
by \eqref{left} coincides 
with the usual one. If $a_i \ne -1$, the map $M_\ba \ni y 
  \mapsto x_i y \in M_{\ba + \varepsilon_i}$ is bijective. 
In fact, its inverse is 
$\frac1{a_i+1} \, \partial_i : M_{\ba + \varepsilon_i} \to M_\ba$.  
If $M$ is a finitely generated left $A$-module, then 
$\dim_k M_\ba < \infty$ for all $\ba \in \ZZ^n$. 
(In fact, if $V \subset M_\ba$ is a $k$-vector subspace 
and $M' := A V \subset M$ is the submodule  
generated by $V$, then $M' \cap M_\ba = V$ by the construction. 
Since $M$ is a noetherian $A$-module, $M_\ba$ is finite dimensional.) 
Hence $M_\rat(-\b1)$ is a straight $S$-module in this case. 

While $M_\rat = 0$ in many cases, we have the following.

\begin{prop}\label{dense}
Let $\Amd$ be the category of finitely generated left $A$-modules. 
Then $(-)_\rat(-\b1) : \Amd \to \Str$ is a dense functor. 
\end{prop}

\begin{proof} 
If $f: M \to N$ is an $A$-homomorphism, we have $f(M_\ba) \subset N_\ba$ 
for all $\ba \in \ZZ^n$. So $(-)_\rat(-\b1)$ gives a functor. Next we prove 
the density. Let $M$  be a $\ZZ^n$-graded $S$-module with $M(-\b1) \in \Str$.  
We will define $\partial_i y$ for $y \in M_\ba$ as follows. 

\medskip

$(*)$ If $a_i \ne 0$, the map $M_{\ba - \varepsilon_i} \ni z \mapsto x_i z 
\in M_{\ba}$ is bijective, hence there is a unique element 
$y' \in M_{\ba - \varepsilon_i}$ such that $x_i y'=y$. 
Set $\partial_i y := a_i y'$.   If $a_i =0$, we set $\partial_i y=0$. 

\medskip

\noindent 
It is easy to check that $(*)$ makes $M$ a left $A$-module with  $M = M_\rat$. 
\end{proof}

In the situation of the proof of Proposition~\ref{dense}, 
$(*)$ is {\it not} a unique way to make $M$ an 
$A$-module. Consider the case $n=1$ (i.e., $S=k[x]$).   
Set $M := A/A \, x\partial$. Then $M$ has a $k$-basis 
$\{ 1, x, x^2, \ldots, \partial, \partial^2, \ldots \}$. 
So $M = M_\rat$ and $M(-\b1) \cong \can \oplus \E(k)$ 
as $S$-modules. Since $\partial M_0 \ne 0$, 
the $A$-module structure of $M$ is not given by $(*)$.

We say a finitely generated left $A$-module $M$ is a straight 
$A$-module if $M = M_\rat$ and its $A$-module structure is given by $(*)$.  
If $M$ and $N$ are straight $A$-modules, then an $A$-homomorphism 
$f : M \to N$ is nothing other than a $\ZZ^n$-graded $S$-homomorphism. 
Thus the category $\Str_A$ of straight $A$-modules is equivalent to $\Str_S$. 

A local cohomology module $H_I^i(S)$ has 
a natural $A$-module structure for any ideal $I$ (cf. \cite{Lyu}). 
In the monomial ideal case, we have the following. 

\begin{prop}\label{H_I as A-module}
Let $I_\Delta$ be a squarefree monomial ideal. 
Then $H_{I_\Delta}^i(S)$ is a straight $A$-module (i.e., the $A$-module 
structure is given by $(*)$). 
\end{prop}

\begin{proof}
Recall that 
$H_{I_\Delta}^i(S)$ is the $i^{\rm th}$ cohomology of the \v{C}ech complex 
$\C$ with respect to  monomial generators of $I_\Delta$. 
Each term of $\C$ is a direct sum of copies of the localizations 
$S_{x^F}$ of $S$ at $\{x^F, x^{2F}, \ldots \}$.  Note that 
$S_{x^F}(-\b1) \cong \Zc(S(-F^\cmpl))$ is a straight $S$-module, and 
its $A$-module structure as a localization of $S$ 
is give by $(*)$. Thus $\C$ is a complex of straight $A$-modules. 
The natural $A$-module structure of $H_{I_\Delta}^i(S) \cong H^i(\C)$ 
is given in this way. So we are done. 
\end{proof}

Usually, the canonical module $\can$, which is a straight $S$-module, 
is regarded as a {\it right} $A$-module using Lie differentials. 
So it seems that a straight $S$-module $M$ itself (not the shifted module 
$M(\b1)$) should be a right $A$-module. 

 For a  right $A$-module $M$, consider an $A$-submodule  
\begin{equation}\label{right}
M_\rat : = \bigoplus_{\ba \in \ZZ^n} M_\ba, 
\quad 
\text{where} \quad  M_\ba = \{ \, y \in M \mid 
\text{$y \, (x_i\partial_i)  =  -a_iy$ for all $i$} \, \}. 
\end{equation}
If $M$ is finitely generated, $M_\rat$ is a straight $S$-module 
by the same argument as left $A$-modules. 
Conversely, any straight $S$-module can be a right $A$-module 
with $M_\rat = M$ as Proposition~\ref{dense} . 
The right $A$-module $\can$ satisfies $\can = (\can)_\rat$, and 
the $\ZZ^n$-grading given by \eqref{right} coincides with 
the one given by $\can \cong S(-\b1)$. 

It is well-known that a left $A$-module $M$ can be viewed as 
a right $A$-module, if we set $yx_i = x_iy$ and 
$y \partial_i = - \partial_i y$ for $y \in M$. 
When we regard $M$ as a right $A$-module in this way, we denote it by $M_A$.   
Then $(M_A)_\rat \cong M_{\rat}(-\b1)$ as $S$-modules. 
So the degree shifting by $-\b1$ also appears here. 
It is also noteworthy that, for a $\ZZ^n$-graded $S$-module $M$,  
the shifted module $M(-\b1)$ is straight if and only if so is 
the graded Matlis dual ${}^*\Hom_S(M, \E(k))$.  
Related arguments for straight modules over a normal semigroup ring 
are found in \S6 of \cite{Y3}. 

\begin{prop}\label{M_rat}
If $M$ is a finitely generated left $A$-module, 
then $M_\rat$ is holonomic. 
\end{prop}

\begin{proof} 
We may assume that $M = M_\rat$. Consider the filtration 
$\Gamma_0 \subset \Gamma_1 \subset \cdots$ with 
$\Gamma_i := \sum_{|\ba| \leq i} M_\ba$ of $M$. 
Then $F_i \Gamma_j \subset \Gamma_{i+j}$ for all $i, j \geq 0$, 
where $\{F_i\}$ is the Bernstein filtration of $A$. 
Hence $\grF M$ has a $\grF A$-module structure. 
Moreover, $\grF M$ is a $\ZZ^{2n}$-graded $\grF A = 
k[\bar{x}_1, \ldots, \bar{x}_n, \bar{\partial}_1, \ldots, 
\bar{\partial}_n]$-module such that the degree of the image of 
$y \in M_{\ba}$, $\ba \in \ZZ^n$, in $\grF M$ is $\bb \in \ZZ^{2n}$, where 
$$b_i = \begin{cases}
a_i & \text{if $i \leq n$ and $a_i \geq 0$,}\\
-a_{i-n} & \text{if $i > n$ and $a_{i-n} < 0$,}\\
0 & \text{otherwise.}
\end{cases}$$ 
It is easy to see that $\grF M$ is a squarefree $\grF A$-module,  
in particular, finitely generated. 
If $[\grF M]_I \ne 0$ for some $I \subset [2n]$, 
then $|I| \leq n$. Thus $\dim_{\grF A} (\grF M) 
\leq n$ (if $M \ne 0$, $\dim_{\grF A} (\grF M) = n$), 
that is, $M$ is holonomic. 
\end{proof}

Let $M$ be a finitely generated left $A$-module.  
Then $M$ admits a good filtration $\{ \Gamma_i \}_{i \geq 0}$, that is, 
the associated graded module $\grF M := 
\bigoplus_{i \geq 0} \Gamma_i / \Gamma_{i-1}$ is a finitely generated 
$\grF A$-module (cf. \cite{Bj}). 
We denote the set of minimal associated primes of $\grF M$ as a 
$\grF A$-module by $\SS(M)$. For $Q \in \SS(M)$, 
we denote the multiplicity of the $\grF A$-module $\grF M$ at $Q$ by $e_Q(M)$  
(cf. \cite[A.3]{Bj1}). It is known that $\SS(M)$ and $e_Q(M)$ do not 
depend on the particular choice of a good filtration of $M$.

For $F \subset [n]$, we denote the monomial prime ideal   
$(\bar{x}_i \mid i \not \in F) + (\bar{\partial}_j \mid j \in F)$
of $\grF A = k[\bar{x}_1, \ldots, \bar{x}_n, 
\bar{\partial}_1, \ldots, \bar{\partial}_n]$ by $Q_F$. 
It is easy to see that $Q_F$ is 
an involutive ideal (i.e., closed under the Poisson product, see 
\cite[A.3]{Bj1}) of dimension $n$. Conversely, every involutive monomial 
prime ideal of dimension $n$ is of the form $Q_F$ for some $F$.

\begin{prop}\label{SS}
Let $M$ be a finitely generated left $A$-module. Then
$$\SS(M_\rat) = \{ Q_F \mid M_{F-\b1} \ne 0 \} \quad \text{and} \quad    
e_{Q_F}(M_\rat) = \dim_k M_{F-\b1}.$$ Here $F$ represents the squarefree 
vector whose support is $F \subset [n]$.   
\end{prop}

We need the following lemma. 

\begin{lem}\label{simple}
Let $M$ be a finitely generated left $A$-module with $M = M_\rat$. 
If $M$ is not a straight $A$-module, then $M$ is not simple as an $A$-module. 
\end{lem}

\begin{proof}
Since $M$ is not straight, there are some $\ba \in \ZZ^n$, $i \in [n]$, and 
$y \in M_\ba$ such that $a_i = 0$ and $\partial_i y \ne 0$. Let 
$N := A (\partial_i y)$ be the submodule of $M$. Since $x_i (\partial_i y) = 
a_i y = 0$, we have $y \not \in N$. Hence $N \ne 0, M$. 
\end{proof}

\noindent{\it Proof of Proposition~\ref{SS}.}
We may assume that $M = M_\rat$. Note that the submodule $N$ constructed 
in the proof of Lemma~\ref{simple} satisfies $N = N_\rat$. (More generally, 
if $M = M_\rat$, any submodule $M'$ of $M$ satisfies $M' = (M')_\rat$.) 
By Lemma~\ref{simple}, we have a filtration 
$0 = M_0 \subset M_1 \subset \cdots \subset M_t = M$ such that 
$M_i = (M_i)_\rat$ and $M_i/M_{i-1}$ is straight for each $i$. 
Recall that $\Str_A \cong \Str_S \cong \Sq_S$ and a simple object in $\Sq$ 
is isomorphic to $(S/P_F)(-F)$ for some $F \subset [n]$. 
So we may assume that $M_i/M_{i-1} \cong \Zc(S/P_F(-F))(\b1) =: L[F]$. 
Take the filtration $\Gamma$ of $L[F]$ given in the proof of 
Proposition~\ref{M_rat}. Then we have $\grF (L[F]) = (\grF A)/Q_F$.  
Hence $\SS(L[F]) =\{ Q_F\}$ and $e_{Q_F}(L[F]) = 1$. 
On the other hand, we have $\dim_k L[F]_{F' -\b1} = \delta_{F, F'}$ 
for all $F' \subset [n]$.  Since $e_{Q_F}(-)$ is additive, we are done. 
\qed 

\bigskip

The characteristic cycle of $H_{I_\Delta}^i(S)$ as an $A$-module 
(i.e.,  $e_{Q_F} (H_{I_\Delta}^i(S))$ for $Q_F \in \SS(H_{I_\Delta}^i(S))$)  
was studied in \cite{Mon}, but we will give another approach here. 
The next corollary shows that Hochster's formula (\cite[Theorem~II.4.1]{St}) 
on the Hilbert function of $H_\m^i(S/I_\Delta)$ is also a formula on the 
characteristic cycle of $H_{I_\Delta}^j(S)$. 

\begin{cor}\label{Hoc}
Let $I_\Delta$ be the Stanley-Reisner ideal of a simplicial complex 
$\Delta \subset 2^{[n]}$. For all $F \subset [n]$ and all $i \geq 0$, we have  
$$e_{Q_F}(H_{I_\Delta}^i(S)) = \dim_k \rH_{n-i+|F|+1}(\lk_\Delta F; k),$$
where $\lk_\Delta F = \{ G \subset [n] \mid \text{$G \cap F = \emptyset$ and 
$F \cup G \in \Delta$}\}$. 
\end{cor}
 
\begin{proof}
By Propositions~\ref{H_I as A-module} and \ref{SS}, 
we have $e_{Q_F}(H_{I_\Delta}^i(S)) = \dim_k [H_{I_\Delta}^i(\can)]_F$.  
But $\dim_k [H_{I_\Delta}^i(\can)]_F = 
\dim_k \rH_{n-i+|F|+1}(\lk_\Delta F; k)$ by Terai's formula (\cite{T}). 
\end{proof}

\begin{rem}
The relation between Hochster's formula and Terai's formula 
is explained by the isomorphisms 
$\Nc(H_{I_\Delta}^i(\can)) \cong \Ext_S^i(S/I_\Delta, \can) \cong 
H_\m^{n-i}(S/I_\Delta)^*$, where $(-)^*$ means the graded Matlis dual.
\end{rem}

\section*{Acknowledgments}
Main parts of this research were done during my stay at 
Mathematical Science Research Institute. I am grateful to this institute 
for warm hospitality.  Special thank are due to Professor David Eisenbud for 
valuable suggestions concerning the Bernstein-Gel'fand-Gel'fand 
correspondence. I also thank Professors Mitsuyasu Hashimoto, Ezra Miller, 
Yuji Yoshino, and the referee for helpful comments.

\end{document}